\documentstyle[11pt]{article}
\newfont{\frak}{eufm10 scaled\magstep1}
\newfont{\bbb}{msbm10 scaled\magstephalf}

\def\thebibliography#1{
\subsection*{R\'ef\'erences bibliographiques\markboth{\sc R\'ef\'erences bibliographiques}{\sc R\'ef\'erences bibliographiques}}
 \list
 {[\arabic{enumi}]}{\settowidth\labelwidth{[#1]}\leftmargin\labelwidth
 \advance\leftmargin\labelsep
 \usecounter{enumi}}
 \def\newblock{\hskip .11em plus .33em minus -.07em}
 \sloppy
 \sfcode`\.=1000\relax}

\newtheorem{thm}{Theorem}
\newtheorem{prop}[thm]{Proposition}
\newtheorem{fthm}{Th\'eor\`eme}
\newtheorem{fdefn}[fthm]{D\'efinition}
\newtheorem{fex}[fthm]{Exemple}

\newcommand{\qed}{\nolinebreak\hfill{$\Box$}\par\vspace{0.5\parskip}}
\newcommand{\qedef}{\nolinebreak\hfill{$\triangle$}\par\vspace{0.5\parskip}}
\newcommand{\qee}{\nolinebreak\hfill{$\Diamond$}\par\vspace{0.5\parskip}}

\def\ut{\tilde{U}}
\def\vt{\tilde{V}}
\def\ft{\tilde{f}}

\def\D{\Delta}
\def\G{\Gamma}
\def\ga{\Gamma_{\alpha}}

\def\ot{\tilde{\omega}}
\def\ota{\ot_{\alpha}}

\def\fia{\phi_{\alpha}}
\def\fib{\phi_{\beta}}
\def\ua{U_{\alpha}}
\def\uta{\ut_{\alpha}}

\def\ub{U_{\beta}}
\def\gab{g_{\alpha\beta}}

\def\C{\mbox{\bbb{C}}}
\def\Q{\mbox{\bbb{Q}}}
\def\R{\mbox{\bbb{R}}}
\def\Z{\mbox{\bbb{Z}}}
\def\rk{\R^k}

\def\zs{|z|^2}
\def\ws{|w|^2}
\def\z1s{|z_1|^2}
\def\ztwos{|z_2|^2}

\def\e{e^{2\pi i\theta}}
\def\est{e^{2\pi i\theta\frac{s}{t}}}
\def\eones{e^{\frac{2\pi i\theta}{s}}} 

\def\d{\mbox{\frak d}}
\def\n{\mbox{\frak n}}
\def\t{\mbox{\frak t}}
\def\ddu{\d^*}
\def\ndu{\n^*}
\def\tdu{\t^*}

\def\X{\mbox{X}}

\begin{document}

\noindent {\large{\sc Sur une
g\'en\'eralisation de la notion de $V$-vari\'et\'e}}\\

\noindent{\bf Elisa PRATO}\\

\noindent \small{Laboratoire Dieudonn\'e, Universit\'e de Nice,
Parc Valrose, 06108 Nice Cedex 2, France\\
Courriel: {\tt elisa@math.unice.fr}\\ Web: {\tt
http://www-math.unice.fr/\~{}elisa/home.html}}\\

\noindent {\bf R\'esum\'e.}  {\small 
Nous consid\'erons un espace
topologique qui est localement isomorphe au quotient de $\rk$ par
l'action d'un groupe discret et nous l'appelons {\em
quasi-vari\'et\'e} de dimension $k$. Les quasi-vari\'et\'es
g\'en\'eralisent les vari\'et\'es et les $V$-vari\'et\'es et repr\'esentent le cadre naturel pour la r\'eduction symplectique par rapport \`a l'action induite d'un sous-groupe de Lie, compact ou non, d'un tore. Nous d\'efinissons les {\em quasi-tores}, les actions hamiltoniennes de quasi-tores et l'application moment sur une quasi-vari\'et\'e symplectique, et nous montrons que tout polytope convexe simple, rationnel ou non, est l'image de l'application moment pour l'action d'un quasi-tore sur une quasi-vari\'et\'e.}\\

\noindent {\sc On a generalization of the
notion of orbifold}\\

\noindent {\bf Abstract.}  {\small We
consider a topological space which is locally isomorphic to the
quotient of $\rk$ by the action of a discrete group and we call it
{\em quasifold} of dimension $k$. Quasifolds generalize manifolds and
orbifolds and represent the natural framework for performing
symplectic reduction with respect to the induced action of any Lie
subgroup, compact or not, of a torus.  We define {\em quasitori},
Hamiltonian actions of quasitori and the moment mapping for symplectic
quasifolds, and we show that every simple convex polytope, rational or
not, is the image of the moment mapping for the action of a quasitorus
on a quasifold.}\\

\subsection*{Introduction}
Soit $M$ une vari\'et\'e symplectique compacte et connexe
et soit $T$ un tore agissant sur $M$ de fa\c{c}on hamiltonienne.
Alors, d'apr\`es le th\'eor\`eme de convexit\'e d'Atiyah et Guillemin-Sternberg
\cite{a,gs}, l'image de l'application moment correspondante
est un polytope convexe {\em rationnel}. Si l'action de $T$ est
effective et $\dim{M}=2\dim{T}$ alors, d'\`apres un th\'eor\`eme de
Delzant \cite{d}, le polytope image est un polytope rationnel simple
satisfaisant \`a une certaine condition d'int\'egralit\'e, et ce
polytope d\'etermine la vari\'et\'e \`a symplectomorphisme
\'equivariant pr\`es.  Au cours de la d\'emonstration de ce
th\'eor\`eme, Delzant donne une construction explicite de la
vari\'et\'e \`a partir du polytope; cette construction se base
sur la technique de r\'eduction symplectique. Les m\'ethodes d'Atiyah,
Guillemin-Sternberg et Delzant ont succesivement \'et\'es adapt\'ees
par Lerman-Tolman \cite{lt} au cas des $V$-vari\'et\'es, mais les
polytopes image, obtenus dans ce cas, sont toujours des polytopes
rationnels. Il est d'ailleurs tr\`es naturel de se demander si tout
polytope simple, {\em rationnel ou non}, est l'image de l'application
moment pour un espace symplectique ad\'equat.  Pour pouvoir r\'epondre
\`a cette question nous consid\'erons des espaces topologiques qui g\'en\'eralisent les vari\'et\'es et les $V$-vari\'et\'es, mais qui ne sont plus forcement des espaces s\'epar\'es.  Ces espaces, que nous appelons {\em quasi-vari\'et\'es}, permettent l'utilisation de la technique de r\'eduction symplectique sous des hypoth\`eses assez faibles; par cons\'equent, ils donnent le moyen d'\'etendre la construction de Delzant au cas des polytopes non rationnels. L'unicit\'e dans le th\'eor\`eme de Delzant n'est plus valable
dans ce contexte (voir exemple~\ref{quasisphere2}). Dans le cas
o\`u $\D$ est rationnel cette construction nous donne une famille de
$V$-vari\'et\'es, en accord avec \cite{lt}; si $\D$ v\'erifie aussi la
condition d'int\'egralit\'e de Delzant, on trouve parmi cette famille
une vraie vari\'et\'e, en accord avec \cite{d}.

Remarquons que l'espace des orbites de l'action d'un groupe discret
sur une vari\'et\'e a \'et\'e etudi\'e par Connes \cite[chapitre
II]{c} dans le cadre de la g\'eom\'etrie non commutative. Les
quasi-vari\'et\'es sont aussi reli\'ees \`a la g\'eom\'etrie des
quasi-cristaux \cite{se}.

Nous renvoyons le lecteur \`a \cite{p1} pour les preuves des r\'esultats \'enonc\'es, ainsi que pour de nombreux exemples; dans une suite \`a \cite{p1} on trouvera un traitement plus detaill\'e des
propri\'et\'es de convexit\'e de l'application moment. Finalement dans
un travail en collaboration avec Battaglia \cite{bp} nous introduisons
des structures complexes et des structures de K\"ahler sur les
quasi-vari\'et\'es, et nous montrons que les espaces du
th\'eor\`eme~\ref{toutsimple} peuvent \^etre consid\'er\'es comme la
g\'en\'eralisation naturelle des vari\'et\'es toriques qu'on associe
\`a tout polytope convexe simple qui est rationnel.

\subsection*{D\'efinitions et r\'esultats}
Commen\c{c}ons par donner la d\'efinition de mod\`ele 
de quasi-vari\'et\'e et de diff\'eomorphisme entre mod\`eles.
 
\begin{fdefn}[Mod\`ele]{\rm
    Soit $\ut$ une vari\'et\'e connexe de
    dimension $k$ et soit $\G$ un groupe discret agissant
    diff\'erentiablement sur $\ut$
    de fa\c{c}on que l'ensemble des points o\`u l'action est libre,
    soit connexe et dense. Un {\em mod\`ele} de dimension $k$, 
    $\ut/\G$, est l'espace des orbites de l'action de $\G$ sur $\ut$,
    muni de la topologie quotient.}\qedef\end{fdefn}
On peut toujours supposer, en passant \'eventuellement au rev\^{e}tement universel, que $\ut$ soit simplement connexe, ce que nous ferons par la suite.
\begin{fdefn}[Diff\'eomorphisme de mod\`eles]{\rm
    Soient $p\,\colon\ut\rightarrow\ut/\G$ et
    $q\,\colon\vt\rightarrow\vt/\D$ deux mod\`eles.
    On dira qu'une application
    $f\,\colon \ut/\G\longrightarrow \vt/\D$ est un
    {\em diff\'eomorphisme de mod\`eles} s'il existe un
    diff\'eomorphisme
    $\ft\,\colon\ut\longrightarrow\vt$ et un isomorphisme
    $F\,\colon\Gamma\rightarrow\Delta$
    tels que $q\circ \ft=f\circ p$ et
    $\ft(\gamma\cdot-)=F(\gamma)\cdot\ft(-)$, $\gamma\in\G$.
    Nous dirons alors que $\ft$ est un {\em relev\'e} de $f$.
}\qedef\end{fdefn}
On montre aisement que deux relev\'es d'un diff\'eomorphisme
$f$ co\"{\i}ncident \`a moins d'une multiplication \`a droite par
un \'el\'ement de $\G$ (ou \`a gauche par un \'el\'ement de $\D$).
Nous pouvons maintenant d\'efinir une quasi-vari\'et\'e
\`a partir des deux d\'efinitions pr\'ec\'edentes.
\begin{fdefn}[Quasi-vari\'et\'e] {\rm Soit $M$ un espace topologique.
    Un {\em atlas de quasi-vari\'et\'e} de dimension $k$ 
     sur $M$ est une collection d'ouverts ${\cal A}= \{\;
    \ua\,|\,\alpha\in A\;\}$, dits {\em cartes}, ayant 
    les propri\'et\'es suivantes:
    \begin{enumerate}
    \item les $\ua$ recouvrent $M$;
    \item pour tout $\alpha\in A$ il existe un mod\`ele
    $\uta/\ga$, o\`u $\uta$ est un sous-ensemble ouvert, 
    connexe et simplement connexe de $\rk$, et un hom\'eomorphisme
    $\fia\,\colon\uta/\ga\longrightarrow\ua$;
    \item pour tous $\alpha, \beta \in A$ tels que
    $\ua\cap\ub\neq\emptyset$ l'application
    $$\gab=\fib^{-1}\circ\fia\,\colon\fia^{-1}(\ua\cap\ub)
    \longrightarrow\fib^{-1}(\ua\cap\ub)$$
    est un diff\'eomorphisme de mod\`eles. Nous dirons alors que
    $\gab$ est un {\em changement de cartes}.
    \end{enumerate}  
    Deux atlas sur $M$ sont \'equivalents si l'union des leurs cartes
    est encore un atlas de $M$. Une {\em structure de
    quasi-vari\'et\'e} sur $M$ est une classe d'\'equivalence
    d'atlas; nous dirons alors qu'un espace $M$ avec une
    structure de quasi-vari\'et\'e est une {\em quasi-vari\'et\'e}.
}\qedef\end{fdefn}
Si les groupes $\ga$ sont finis on retrouve ici la d\'efinition classique de $V$-vari\'et\'e, s'ils sont triviaux on retrouve la d\'efinition de vraie vari\'et\'e. 
\begin{fex}[Quasi-sph\`ere]\label{quasisphere}{\rm 
    Soient $s,t$ deux nombres r\'eels positifs tels que $s/t\notin\Q$.
    Consid\'erons $\C^2$ avec sa forme symplectique canonique
    $\omega_0=\frac{1}{2\pi i}(dz_1\wedge d\bar{z}_1 +
    dz_2\wedge d\bar{z}_2)$ et l'action de $\R$: 
    $(\theta, (z_1,z_2)) = ( \e z_1, \est z_2)$ d'application moment
    $\Psi(z_1,z_2)=\z1s +\frac{s}{t}\ztwos-s$.
    Consid\'erons le niveau r\'egulier $\Psi ^{-1} (0)$; alors
    l'espace des orbites $M=\Psi^{-1}(0)/\R$ est une quasi-vari\'et\'e
    de dimension $2$.
}\qee\end{fex}
Soit $M$ une quasi-vari\'et\'e.
Pour d\'efinir une g\'eom\'etrie sur $M$
nous proc\'edons de la fa\c{c}on suivante:
d'une mani\`ere g\'en\'erale, un objet g\'eom\'etrique sur une
quasi-vari\'et\'e est la donn\'ee d'un objet g\'eom\'etrique
sur tout ouvert $\uta$, invariant par rapport \`a l'action du groupe discret $\ga$ et satisfaisant \`a des conditions de compatibilit\'e. 
Par exemple, une forme
diff\'erentielle de degr\'e $h$, $\omega$, sur $M$ est la donn\'ee,
pour tout $\alpha\in{\cal A}$, d'une forme diff\'erentielle 
$\ga$-invariante de degr\'e $h$, $\ota$, sur $\uta$ qui
se transforme de mani\`ere compatible aux changements de cartes
$\gab$; une $2$-forme est dite symplectique
si toute $\ota$ est symplectique (ferm\'ee et non-d\'eg\'en\'er\'ee).
Nous renvoyons le lecteur \`a l'article \cite{p1} pour
les d\'efinitions de: champ de vecteurs, application diff\'erentiable,
diff\'eomorphisme, image r\'eciproque d'une forme, 
image directe d'un champ de vecteurs, 
diff\'erentiel et produit int\'erieur.
L'analogue du tore dans la g\'eom\'etrie des
quasi-vari\'et\'es s'appelle quasi-tore.
Soit $\d$ un espace vectoriel de dimension $n$.
\begin{fdefn}[Quasi-r\'eseau, quasi-tore, quasi-alg\`ebre de Lie]{\rm
  Un {\em quasi-r\'eseau}, $Q$, dans $\d$ est le $\Z$-module
  engendr\'e par un ensemble de vecteurs
  $X_1,\ldots,X_d$ qui engendrent $\d$. Nous appelons {\em quasi-tore}   de dimension $n$ le groupe quotient $D=\d/Q$.
  La {\em quasi-alg\`ebre} de Lie du quasi-tore $D$
  est l'espace vectoriel $\d$.}
\qedef\end{fdefn}
Un quasi-tore est evidemment une quasi-vari\'et\'e et les operations de groupe sont diff\'erentiables. Si $Q$ est un vrai r\'eseau (par exemple si $d=n$) on retrouve dans la d\'efinition pr\'ec\'edente le tore et son alg\`ebre de Lie. Les quasi-tores de dimension un ont
\'et\'es \'etudi\'es par Donato, Iglesias et Lachaud
\cite{di,i, il}; Iglesias introduit \`a cette
occasion la terminologie {\em tores irrationnels}.
\begin{fex}[Quasi-cercle]\label{quasicercle}{\rm Consid\'erons
dans $\R$ le quasi-r\'eseau $Q=s\Z+t\Z$, $s/t\notin\Q$. Alors
$D^1=\d/Q$ est un quasi-tore de dimension $1$. \qee}\end{fex}
Les quasi-tores apparaissent de fa\c{c}on tr\`es naturelle,
comme le montre la proposition suivante.
\begin{prop}
Soit $T$ un tore et soit $N$ un sous-groupe de Lie de $T$.
Alors $T/N$ est un quasi-tore de dimension $n=\dim{T}-\dim{N}$.
\qed\end{prop}
Par exemple le quotient du tore de dimension $2$
par une droite de pente irrationnelle $s/t\notin\Q$ est
le quasi-tore de l'exemple~\ref{quasicercle}.
\begin{fdefn}[Action diff\'erentiable]{\rm
    Une {\em action diff\'erentiable} d'un quasi-tore $D$ sur 
    une quasi-vari\'et\'e $M$ est une application
    diff\'erentiable $\tau\,\colon\, D\times M \longrightarrow M$
    telle que $\tau(d_1\cdot d_2,m)=\tau(d_1,\tau(d_2,m))$, et
    $\tau(1_{D},m)=m$ pout tout $d_1, d_2 \in D$ et $m\in M$.}
\qedef\end{fdefn}
Etant donn\'e une action diff\'erentiable d'un quasi-tore sur une
quasi-vari\'et\'e il est toujours possible d'associer \`a tout $X\in\d$, comme dans le cas des vari\'et\'es, 
un champ de vecteurs sur $M$, dit
{\em champ fondamental de l'action}, not\'e $\X_M$.
\begin{fdefn}[Action hamiltonienne, application moment]{\rm
    Une action diff\'erentiable d'un quasi-tore $D$ sur 
    une quasi-vari\'et\'e $M$ est dite
    {\em hamiltonienne} si elle
    respecte la forme symplectique et s'il existe une
    application diff\'erentiable 
    $D$-invariante $\Phi\,\colon\, M\longrightarrow\ddu$,
    dite {\em application moment}, ayant la propri\'et\'e que
    $\imath(\X_M)\omega=d<\Phi,X>$, pour tout $X\in\d$.
\qedef}\end{fdefn}
\begin{fex}\label{quasisphere2}{\rm 
  Consid\'erons la quasi-sph\`ere de l'exemple~\ref{quasisphere}
  et le quasi-cercle de l'exemple~\ref{quasicercle}.
  L'application $\tau([\theta],[z:w])=[\eones z:w]$
  d\'efinit une action hamiltonienne de $D^1$ sur $M$,
  d'application moment $\Phi([z:w])=\frac{\zs}{s}=1-\frac{\ws}{t}$.
  Remarquons que $\Phi(M)=[0,1]$ (indep\'endamment de $s$ et $t$)
  comme pour la rotation de la sph\`ere de $\R^3$ autour
  de l'axe $Oz$.}\qee\end{fex}
Pour obtenir beaucoup d'autres exemples
d'actions hamiltoniennes de quasi-tores sur des quasi-vari\'et\'es
il suffit de prendre la r\'eduction symplectique par rapport \`a
un sous-groupe de Lie d'un vrai tore agissant sur une
vraie vari\'et\'e symplectique.
\begin{fthm}[R\'eduction]
  Soit $T$ un tore d'alg\`ebre de Lie
  $\t$, soit $T\times X \longrightarrow X$ une action hamiltonienne de 
  $T$ sur une vari\'et\'e symplectique $X$ et supposons que l'application
  moment $J\,\colon X \longrightarrow\tdu$ soit propre. 
  Consid\'erons l'action induite d'un sous-groupe de Lie $N$ de
  $T$ et supposons que $0$ soit une valeur r\'eguli\`ere   de l'application
  moment correspondante: $\psi=i^*\circ J\,\colon X \longrightarrow \ndu$ ($\n$
  d\'enote l'alg\`ebre de Lie de $N$ et $i$ l'inclusion dans
  $\t$.) Alors $M=\psi^{-1}(0)/N$ est une quasi-vari\'et\'e symplectique
  de dimension $\dim{X}-2\dim{N}$ et l'action induite de $T/N$ sur $M$ 
  est hamiltonienne.
\qed \end{fthm}
Ce dernier th\'eor\`eme permet d'\'etendre la construction
de Delzant et donc de montrer que tout polytope simple $\D\subset\ddu$
est l'image par l'application moment d'une quasi-vari\'et\'e.
Rappelons qu'un polytope convexe est dit
simple si de chaque sommet sont issues
exactement $n$ ar\^etes (ici nous supposons
par simplicit\'e $\dim{\D}=n$). Nous avons alors:
\begin{fthm}\label{toutsimple}
Pour tout polytope convexe simple $\D\subset\ddu$ il existe un
quasi-tore de dimension $n$ et de quasi-alg\`ebre de Lie $\d$, $D$,
une quasi-vari\'et\'e symplectique compacte de dimension $2n$,
$M$, et une action hamiltonienne effective de $D$ sur $M$
telle que l'image de l'application moment correspondante soit $\D$.
\qed\end{fthm}

\end{document}